\documentstyle[11pt]{article}
\title{Finding Large Monochromatic Diameter Two Subgraphs}
\author{Tom Fowler \thanks{Research supported in part by ONR under
contract N00014-93-1-0325}\\ School of Mathematics  \\ Georgia Institute of
Technology \\ Atlanta, Georgia 30332}

\newtheorem{lemma}{Lemma}
\newtheorem{theorem}{Theorem}

\begin{document}
\maketitle





\begin{abstract}
\noindent
Given a coloring of the edges of the complete graph on $n$ vertices in
$k$ colors, by considering the neighbors of an arbitrary vertex it
follows that there is a monochromatic diameter two subgraph on at least
$1+(n-1)/k$ vertices. We show that for $k \ge 3$ this is asymptotically
best possible, and that for $k=2$ there is always a monochromatic
diameter two subgraph on at least $\lceil {3 \over 4} n \rceil$ vertices, which again, is best
possible.
\end{abstract}

\section{Introduction}

\hspace{.2 in} Ramsey's Theorem implies that for every positive
integer $l$, there is an integer $n$ such that for every two-coloring of
the edges of a complete graph on at least $n$ vertices, there is a
monochromatic clique (or equivalently, a monochromatic diameter one subgraph)
containing at least $l$ vertices. Many variations of this
classical problem have been considered, and there are several books
devoted exclusively to the subject of Ramsey Theory (\cite{GRS,NR}).
I consider
another variation that was proposed by Paul Erd\H os \cite{Erd}, which apparently has not been studied previously: Given
a complete graph $K\sb n$ on $n$ vertices whose edges are colored red or blue, what is the largest (in terms of number of vertices) monochromatic diameter two
subgraph that $K\sb n$ is guaranteed to contain?
Here monochromatic means that all of the edges have the same color and diameter two means that for every pair of vertices $x,y$ in the subgraph there is a path in the subgraph joining $x$ and $y$ and containing at most two edges.

	More generally, given a complete graph $K\sb n$ on $n$ vertices whose edges are
partitioned into $k$ different sets which we will hereafter refer to as color
classes, what is the largest monochromatic diameter two subgraph that
$K\sb n$ is guaranteed to have. Erd\H os \cite{Erd} noticed that it is always possible
to find a monochromatic diameter two subgraph of size $1+{(n-1)\over k}$ by considering  
at any vertex $v$ the largest monochromatic star centered at
$v$. 
Is it possible to do better? We will show that the answer to this
question is yes if and only if $k=2$.

In section two of this paper we show that
when $k=2$, any edge two-coloring of $K\sb n$ has a
monochromatic diameter two subgraph containing at least
$\lceil {3 \over 4} n \rceil$ vertices but there is a
two-coloring of the edges of $K\sb n$ such that the largest monochromatic
diameter two subgraph contains at most $\lceil {3 \over 4} n \rceil$ vertices.
In section $3$, for $k=3$ colors and for an integer $s$
divisible by $3$, we construct three-colorings of the edges
of $K\sb n$ on $n=2sk+1$ vertices whose largest monochromatic diameter two subgraph has size
$2s+3$. In section $4$, for any
$k \ge 4$ and any positive integer $s$, we induce a partition of the edges of the complete graph on $n=2sk+1$ vertices into $k$ color classes from a construction of $k$ circle graphs, each of which 
has no diameter two subgraph containing more than
$2s+7$ vertices. Below we provide upper bounds when $n$ does
not have the form $2ks+1$.

	When $2ks+1 < n < 2k(s+1)+1$, say $n=2ks+l$, the following construction along with the constructions above show that there is a $k$
coloring of the edges of $K\sb n$ whose largest
monochromatic diameter two subgraph contains at most $2s+6+l$
vertices.
Given an edge $k$-colored $K\sb {n-1}$ on vertices $\{ 1,..,n-1\}$, create an edge $k$-coloring on $K\sb n$ by  
introducing a new vertex labeled $n$, coloring the edge $\{ n-1,n\}$ with any
one of the $k$ colors, and for $1 \le x \le n-2$, coloring the
edge $\{ x,n\}$ the same color as the edge $\{ x,n-1\}$. Then the
largest monochromatic diameter two subgraph in $K\sb n$ contains at most one more vertex
than the largest monochromatic diameter two subgraph in $K\sb {n-1}$ since the vertex $n$ replicates the vertex $n-1$. 
When $k=3$ colors and $6s+1 < n < 6(s+3)+1$ this gives an 
upper bound of $2s+9$.

\section{Two Colors} 

\hspace{.2 in} Given a graph $G=(V,E)$, and a subset $A \subset V$, the graph
$G-A$, will denote the subgraph of $G$ with vertex set $V-A$ and
edge set consisting of every edge in $E$ which has both endpoints in
$V-A$.          
 
First we exhibit edge colorings of $K\sb n$ that have no monochromatic
diameter two subgraphs on more than $\lceil {3 \over 4} n \rceil$ vertices.
Partition $V(K\sb n)$
into $4$ sets $R\sb 1$, $R\sb 2$, $B\sb 1$, $B\sb 2$, each of size $\lfloor
{n \over 4} \rfloor$ or $\lceil {n \over 4} \rceil$. 
Color an edge red if it has one endpoint in $R\sb 1$ and the
other in $B\sb 1$ or if it has one endpoint in $B\sb 1$ and the
other in $B\sb 2$ or one endpoint in $B\sb 2$ and the other
in $R\sb 2$; otherwise color it blue. It is easily verified that
the shortest path with all red edges between any vertex in
$R\sb 1$ and any vertex in $R\sb 2$ has three edges. Similarly the shortest path with
all blue edges between any vertex in $B\sb 1$ and any vertex in $B\sb 2$ has three
edges.  Therefore any monochromatic diameter two subgraph $H$ must be disjoint from one of
$R\sb 1$, $R\sb 2$, $B\sb 1$, and $B\sb 2$, and hence contains at most $\lceil 
{3 \over 4} n \rceil$ vertices.
																		 
	Now we show that any edge two-coloring of $K\sb n$ will
render a monochromatic diameter two subgraph of order at
least $\lceil {3n \over 4} \rceil$. Let $R \bigcup B$ be a partition
of the edges of a graph $G$ into red edges, (those in $R$) and blue edges, (those in $B$).  A path $P$ in $G$ is said to be a red (blue)
path if all of $P$'s edges are red (blue).  A vertex $x\in V(G)$ is said to be a red (blue)
violated vertex in $G$ if there is a $y \in V(G)$ such that there is no red (blue) $x-y$ path in $G$
with two or fewer edges.  Similarly, a pair ${x,y}$ of vertices is said to be a red
(blue) violated pair in $G$ if $ x \not= y $ and there is no red (blue) path in  $G$ joining
$x$ and $y$ containing two or fewer edges.

	\begin{lemma}   Let the edges of $K\sb n$ be colored arbitrarily red or
blue. A vertex $x\in V(K\sb n)$ cannot be both red violated in $K\sb n$
and blue violated in $K\sb n$. Moreover if $w,x$ is a red (blue)  violated pair in
$K\sb n$ and $y$ is a blue (red) violated vertex in $K\sb n$ then the edges $\{ w,y \} $ and
$\{ x,y \}$ have different colors.  \end{lemma}

\underline{Proof.}    First, suppose to the contrary that there are
vertices $r,b$ such that $x,r$ is a red violated pair in $K\sb n$ and $x,b$ is a blue violated pair in $K\sb n$. It follows that $\{ x,r \} \in B$ and $ \{ x,b \} \in R$.
If  $\{ r,b \} \in R$  then $xbr$ is a red path
containing two edges joining $x$ and $r$ contradicting the assumption
that $x,r$ is a red violated pair in $K\sb n$. A similar
contradiction arises if $\{ x,b \} \in B$.
Now suppose that $w,x$ is a red violated pair in $K\sb n$ and that there exists
a vertex $z$ such that $y,z$ is a blue violated pair in $K\sb n$.
Clearly, the edges $\{ w,y \}$ and $\{ x,y \} $ can't both be red.
Now suppose that they are both blue. This forces both of the edges $\{ w,z\} $ and
$\{ x,z\} $ to be red since $y,z$ is a blue violated pair in $K\sb n$ which
contradicts the assumption that $w,x$ is a red violated pair. By
symmetry, the same proof holds if the colors red and blue are
interchanged. \hspace{3.55 in} $\Box$

Lemma $1$ imposes a great deal of structure on the edges joining
vertices in red and blue violated pairs. The key idea of the following
proof is to use this structure to partition the red and blue violated
vertices in an economical way.

\begin{theorem} Let $R\bigcup B$ be an arbitrary partition of the
edges of $K\sb n$. Then there exists a monochromatic diameter two subgraph $H$
with $|V(H)| \ge \lceil {3n \over 4} \rceil$. Moreover, there exists edge two-colorings for which equality is attained. 
\end{theorem}

\underline{Proof.}    If there exists a spanning monochromatic diameter two subgraph $H$,
then we are done. Therefore, we may assume that there is a red violated pair $ r\sb 1,r\sb 2$  in $K\sb n$ and a blue violated pair $ b\sb 1,b\sb 2$ in $K\sb n$. By
Lemma $1$, $ \{ r\sb 1,r\sb 2\} \bigcap \{ b\sb 1,b\sb 2\} = \emptyset$.
By Lemma $1$ we may assume that
$\{ r\sb 1,b\sb 1 \} , \{ r\sb 2,b\sb 2\}\in R$ and $\{ r\sb 1,b\sb
2\} , \{ r\sb 2,b\sb 1 \} \in B$. For $\delta = 1,2$ and $i=1,2,\ldots$
define $R\sb {\delta ,i} $ to be the set of all $x$ that are red violated
vertices in $K\sb n-$$(R\sb {\delta ,1} \bigcup \ldots \bigcup R\sb{\delta ,i-1}) $
and for which  the edge $\{ x,b\sb {\delta} \}$ is red and
$B\sb {\delta ,i} $ to be the set of all $y$ that are blue violated vertices in
$K\sb n-$$(B\sb {\delta ,1} \bigcup \ldots \bigcup B\sb{\delta ,i-1}) $ and for
which the edge $\{ y,r\sb {\delta} \} $ is red.
For $i=1$ interpret the union to be empty and note that $r\sb {\delta} \in
R\sb {\delta,1}$ and $b\sb {\delta} \in B\sb {\delta,1}$ for
$\delta = 1,2 $. Since the number of vertices of $K\sb n$ is
finite there is a largest integer $l$ such that at least one of $R\sb
{1,l}$, $R\sb {2,l}$,$B\sb {1,l}$,$B\sb {2,l}$ is nonempty.

	Refer to the sets $R\sb {1,1}$,$\ldots$,$R\sb {2,l}$ as red sets and the sets $B\sb {1,1}$,\ldots,$B\sb {2,l}$
as blue sets.

\vspace{.15 in}
($1$)	For every $\delta \in \{ 1,2 \} $, $r\sb {\delta} $   ($b\sb {\delta}$) is not in any blue set (red set).
\vspace{.15 in}

\underline{Proof of ($1$).}
We will prove by induction on
$i$ that $r\sb {\delta} \notin B \sb {1,i} \bigcup B\sb
{2,i} $ for every $\delta \in \{ 1,2 \}$ and for every positive
integer $i$. Suppose that $i=1$. Then $r \sb
{\delta}$ in either $B\sb {1,1}$ or in $B\sb {2,1}$
implies that $r\sb {\delta}$ is a blue violated vertex in $K\sb n$.
This however is impossible by Lemma $1$ since $r\sb {\delta}$ in $R\sb {\delta,1}$ implies that
$r\sb {\delta}$ is a red violated vertex in $K\sb n$. Now suppose
that $i \ge 2$. By the induction hypothesis, $\{ r\sb 1,r\sb 2 \}$
is disjoint from $B\sb {\alpha,1} \bigcup \ldots \bigcup B\sb {\alpha,i-1} $ for every
$\alpha \in \{ 1,2 \}$. Therefore, $r\sb {\delta}$ is a red
violated vertex (with $r\sb {3-\delta} $) in $K\sb n-$$(B\sb {\alpha,1}
\bigcup \ldots \bigcup B\sb {\alpha,i-1})$ and hence cannot also
be a blue violated vertex in $K\sb n-(B\sb {\alpha,1} \bigcup \ldots \bigcup
B\sb {\alpha,i-1 })$.
	It follows that $r\sb {\delta} $ is not in $B\sb {\alpha,i}$ for
any $\alpha \in \{ 1,2 \}$. Similar reasoning applies to $b\sb {\delta}$. This
proves ($1$).
	
	    The next proposition will be crucial in guaranteeing the existence of
a sub-collection of the red and blue sets whose union is ``small''
and whose removal results in a monochromatic diameter two subgraph.
\vspace{.1 in}

($2$)	The sets $R\sb {1,1},R\sb {2,1} ,R\sb {1,2}
,\ldots,R\sb{2,l}$,$B\sb {1,1},B\sb {2,1} ,B\sb {1,2},\ldots,B\sb
{2,l}$ are pairwise disjoint.
\vspace{.1 in}

\underline{Proof of ($2$).} We must show that a)  $R\sb{\delta,i} $, $R\sb{\alpha,j} $ are
disjoint unless $\delta = \alpha $ and $i=j$, b) $B\sb{\delta,i}$, $B\sb{\alpha,j}$ are disjoint unless $\delta = \alpha$ and
$i=j$ and c) that any red set and any blue set are disjoint. To
show a) suppose $\delta \not= \alpha $ and $x \in R\sb{\delta,i}$.
Clearly $\{x,b\sb{\delta} \}$ is a red edge. Proposition ($1$) shows that
$b\sb 1,b\sb 2 $ is a blue violated pair in $K\sb n$$-(R\sb{\delta,1} \bigcup
\ldots \bigcup R\sb {\delta,i-1}) $ which by Lemma $1$ implies that $\{ x,b\sb{\alpha} \} $ is a blue edge which
means that $x$ cannot be contained in $R\sb {\alpha,j}$ for any $j$. So
suppose $\delta = \alpha$ and $i<j$.  In this case the definition of $R\sb
{\delta,j}$ shows that $R\sb{\delta,i} $ and $R\sb{\delta,j}$ are disjoint. A similar
argument shows that b) is true. In particular, for every
$\delta \in \{ 1,2 \} $ the only red (blue) set that $r\sb
{\delta}$ ($b\sb {\delta}$) is in is $R\sb {\delta,1}$ ($B\sb
{\delta,1}$) . To show that any red set and any blue
set are disjoint we prove:
\vspace{.15 in}

($3$) For every $\alpha \in \{ 1,2\}$ and for every positive integer $i$ the following two propositions are false:
\vspace{.15 in}

	 A) There is a vertex $w \in R\sb{\alpha,i}$ and a vertex $x$ in the
union of the blue sets such that $w,x$ is a red violated pair in
$K\sb n-(R\sb {\alpha,1} \bigcup \ldots \bigcup R\sb {\alpha,i-1})$.

	 B) There is a vertex $w \in B\sb{\alpha,i}$ and a vertex $x$ in the
union of the red sets such that $w,x$ is a blue violated pair in $K\sb n- (B\sb
{\alpha,1} \bigcup \ldots \bigcup B\sb {\alpha,i-1})$.
	
\vspace{.15 in}

\underline{Proof of ($3$).} Suppose the assertion is false and let $i$ be the smallest
integer such that either A) or B) is true. Without loss of generality suppose that $w$
is in $R\sb {\alpha,i}$, $x$ is in $B\sb {\delta,j} $, $w,x$ is a red
violated pair in $K\sb n-(R\sb {\alpha,1} \bigcup \ldots \bigcup R\sb
{\alpha,i-1})$ and there is a vertex $y$ such that $x,y$ is a blue
violated pair in $K\sb n-(B\sb{\delta,1} \bigcup \ldots \bigcup B\sb {\delta,j-1})$.

	 Claim: $y$ is not in $R\sb {\alpha,1}
\bigcup \ldots \bigcup R\sb{\alpha,i-1} $. If it is then
there is an integer $s<i$ and a vertex $z$ such that $y \in R\sb{\alpha,s}$
and $y,z$ is a red violated pair in $K\sb n-(R\sb{\alpha ,1} \bigcup \ldots \bigcup R\sb{\alpha,s-1}) $.
By the choice of $i$, $z$ is not in  $ B\sb{\delta, 1} \bigcup \ldots \bigcup B\sb{\delta,j-1}$.
 Therefore, $y,z$ is a red violated pair in the graph
		  \noindent $K\sb n-(R\sb {\alpha,1}\bigcup \ldots \bigcup R\sb
{\alpha,s-1} \bigcup B\sb{\delta,1}\bigcup \ldots \bigcup B\sb{\delta,j-1})$ .
Also, $x$ is not in $R\sb{\alpha,1}\bigcup \ldots \bigcup
R\sb{\alpha,i-1}$ because $w,x$ is a red violated pair in
$K\sb n-(R\sb{\alpha,1}\bigcup \ldots \bigcup R\sb{\alpha,i-1})$.
Therefore $x,y$ is a blue violated pair in
$K\sb n-(R\sb{\alpha,1}\bigcup \ldots \bigcup R\sb{\alpha,s-1}\bigcup
B\sb{\delta,1} \bigcup B\sb{\delta,j-1}) $ so $y$
is both red and blue violated in $K\sb n-(R\sb{\alpha,1}\bigcup \ldots$$\bigcup
R\sb{\alpha,s-1}\bigcup B\sb{\delta,1}\bigcup \ldots \bigcup
B\sb{\delta,j-1})$ which contradicts Lemma $1$.
This proves the claim that $y$ is not in $R\sb{\alpha,1}\bigcup \ldots \bigcup
R\sb{\alpha,i-1}$.

If $w$ is not in $B\sb{\delta,1} \bigcup \ldots \bigcup
B\sb{\delta,j-1}$ then $x$ is both red violated (with $w$) and blue violated (with $y$) in
$K\sb n-(R\sb{\alpha,1}\bigcup \ldots \bigcup R\sb{\alpha,i-1}$ $\bigcup B\sb{\delta,1}
\bigcup \ldots \bigcup B\sb{\delta,j-1}) $ which contradicts Lemma $1$.
So we may assume that there is a positive integer $t$ with $t<j$
and a vertex $u$ such that $w \in B\sb{\delta,t} $
and $w,u$ is a blue violated pair in $K\sb n-(B\sb{\delta,1}\bigcup \ldots \bigcup
B\sb{\delta,t-1})$.
Remembering that $x\in B\sb{\delta,j} $ we have that the edges
$\{ w,r\sb{\delta} \}$,$\{ x,r\sb{\delta} \}$ are both colored
red. Since $w,x$ is a red violated pair in $K\sb n-$$(R\sb {\alpha,1}
\bigcup \ldots \bigcup R\sb{\alpha,i-1})$ it must be that $r\sb{\delta} \in R\sb{\alpha,1}
\bigcup \ldots \bigcup R\sb{\alpha,i-1}$ which implies that $\delta=\alpha
$ since the only red set that $r\sb{\delta}$ is in is $R\sb{\delta,1}$.
Claim: $u$ is not in $ R\sb{\alpha,1}\bigcup \ldots \bigcup
R\sb{\alpha,i-1}$.  Proof: Suppose by way of contradiction that the claim is false.
Proposition ($1$) guarantees that $b\sb 1 $,$b\sb 2 $ is a blue violated
pair in the graph $K\sb n-$$(R\sb{\alpha,1} \bigcup \ldots \bigcup R\sb{\alpha,m-1})$
for every positive integer $m$ and allows us to apply the hypothesis
of Lemma $1$ to $u,b\sb 1 ,b\sb 2$  and to $w,b\sb 1 ,b\sb 2$.
Hence, $u$ in $R\sb{\alpha,1} \bigcup \ldots \bigcup R\sb{\alpha,i-1} $
implies that $\{ u,b\sb {3-\alpha} \}$ is a blue edge. Also, $w\in
R\sb{\alpha,i} $ implies that $\{ w,b\sb{3-\alpha} \}$ is a blue edge.
Thus, $wb\sb {3-\alpha} u$ is a blue path containing two edges. Since $w$,$u$ is a blue violated pair in
$K\sb n-$$(B\sb {\delta,1} \bigcup  \ldots \bigcup B\sb{\delta,t-1})$ it must be
that $b\sb {3-\alpha}$ is in $B\sb {\delta,1} \bigcup \ldots \bigcup B\sb {\delta,t-1}$.
 This implies that $3-\alpha = \delta$ because the only
blue set that $b\sb {3-\alpha}$ is in is $B\sb {3-\alpha,1}$.
This contradicts the fact that $\delta = \alpha $ and proves the
claim that $u$ is not in $R\sb {\alpha,1} \bigcup \ldots \bigcup R\sb {\alpha,i-1}$.

	  Also, $x$ is not in $ B\sb{\delta,1} \bigcup \ldots \bigcup B\sb{\delta,t-1} $ because $t<j$
and $x \in B\sb{\delta,j} $. Hence $w$ is both red violated (with $x$) and blue violated (with $u$) in the graph
$K\sb n-(R\sb{\alpha,1}\bigcup \ldots \bigcup R\sb{\alpha,i-1}\bigcup B\sb{\delta,1} \bigcup \ldots \bigcup
B\sb{\delta,t-1})$ and this contradicts Lemma $1$. Thus ($3$) holds.

Now we will show that any red set and any blue set are disjoint.
Suppose to the contrary that $w \in R\sb{\alpha,i}\bigcap B\sb{\delta,j}
$ for some $\alpha , \delta \in \{ 0,1 \} $ and for some positive integers $i,j$. 
Let $x,y$ be vertices such that $w,x$ is a red violated pair in
$K\sb n-(R\sb{\alpha,1}\bigcup \ldots \bigcup R\sb{\alpha,i-1}) $ and $w,y$ is
a blue violated pair in $K\sb n-(B\sb{\delta,1}\bigcup \ldots \bigcup B\sb{\delta,j-1})$.
By ($3$), $x \notin B\sb{\delta,1}\bigcup  \ldots \bigcup B\sb{\delta,j-1},$ $y \notin
R\sb{\alpha,1}\bigcup \ldots \bigcup R\sb{\alpha,i-1}$ and consequently $w$ is both red and blue violated in 
$K\sb n-(R\sb{\alpha,1}\bigcup \ldots \bigcup R\sb{\alpha,i-1}\bigcup
		   B\sb{\delta,1} \ldots \bigcup B\sb{\delta,j-1}) $ which contradicts Lemma $1$. This proves ($2$).

	     From ($2$) it follows that at least one of the following
four propositions holds :

(i) $|R\sb {1,1} \bigcup R\sb {1,2} \bigcup \ldots \bigcup R\sb{1,l} | \le \lfloor {n \over 4} \rfloor$.

(ii) $|R\sb {2,1} \bigcup R\sb {2,2}\bigcup \ldots \bigcup R\sb {2,l} | \le \lfloor {n \over 4} \rfloor$.

(iii) $|B\sb {1,1} \bigcup B\sb {1,2}\bigcup \ldots \bigcup B\sb{1,l} |  \le \lfloor {n \over 4} \rfloor$.
   
(iv) $|B\sb{2,1} \bigcup B\sb {2,2}\bigcup \ldots \bigcup B\sb {2,l} |  \le \lfloor {n \over 4} \rfloor$.

\noindent Without loss of generality assume proposition (i) holds. Let
$H'$ be the subgraph induced by the vertex set
$V(H')$$= V(K\sb n)-$$(R\sb {1,1}\bigcup \ldots \bigcup R\sb {1,l})$.
Claim: For any two vertices $x$ and $y$ in $V(H')$, there exists a red path in
$H'$ connecting 
$x$ to $y$ containing not more than two edges. Proof: Suppose by way of contradiction that $x,y$ is a red
violated pair in $H'$. Since $b\sb 1$,$b\sb 2$ is a blue violated
pair in $H'$, we can conclude from Lemma $1$ (and if necessary, by relabeling)
that the edge joining $x$ to $b\sb 1$ is red.
Now let an integer $s$ be the smallest
integer such that $x,y$ is a red violated pair in $K\sb n-(R\sb{1,1}\bigcup
\ldots \bigcup
R\sb {1,s})$. Clearly such a $s$ exists because $x,y$ is a
red violated pair in $K\sb n-(R\sb{1,1}\bigcup \ldots \bigcup R\sb{1,l} )$.
Now $x \in R\sb {1,s+1} $ by definition of $R\sb {1,s+1}$.
Since $R\sb {1,l+1} = \emptyset $ we must have $s<l$ which implies
that $x \in R\sb {1,1} \bigcup \ldots \bigcup R\sb {1,s+1}
\subseteq R\sb {1,1} \bigcup  \ldots \bigcup R\sb {1,l}$ and this contradicts $x\in V(H')$. This proves the claim that for any two vertices $x$ and $y$ in $H'$ there is a red path in $H'$ with
not more than two edges joining $x$ to $y$.
Define $H=(V(H'),E(H')\bigcap R) $. By this claim and the truth of
proposition (i), $H$ is a monochromatic diameter two subgraph with $|V(H)| \ge
\lceil {3 \over 4} n \rceil$.
	
	The coloring of the edges of $K\sb n$ described at the beginning of this paper shows that $\lceil {3 \over 4} n \rceil$ is the best possible bound. This proves Theorem $1$. Q.E.D.

\section{Three Colors}

\hspace{.2 in}	Let $\{ 0,1,\ldots,n-1\}$ be a ground set. Addition will be modulo $n$
unless otherwise stated. For a positive integer $p$, and integers $a \le b$
with $a,b \in \{ 0,\ldots,n-1 \}$ such that $a \equiv b$$ (\bmod \ p$) we will let $[a,b]\sb p$
denote the set of integers $x$ in $\{ 0,\ldots,n-1 \}$ such that 
$a \le x \le b$, and $x \equiv a \ (\bmod\ p)$.
The notation $[a,b]$ will mean the set $[a,b]\sb 1$ unless
otherwise stated. For $a > b$ let $[a,b] \sb p$ be the set 
$[a,c]\sb p \bigcup [d,b]\sb p$ where $c \le n-1$ is the largest integer
 such that $a \equiv c (\bmod\ p)$ and $d \ge 0$ is the least integer such that 
$b \equiv d (\bmod \ n)$.
If $p=1$ the set $[a,b]$ will also be referred to as an $\it interval$. The quantity $|[a,b]|$, sometimes referred to as the length of the interval $[a,b]$, is
defined to be $1+b-a$ for $a \le b$ and $|[a,n-1]|+|[0,b]|$ otherwise.
Given two sets $S,T
\subset [0,n-1]$, $S\sp c$ will denote the set $[0,n-1]-S$, $S+T$
will denote the set of all elements of the form $s+t$ where $s \in S$
and $t \in T$ and $S-T$ will denote the set of all elements of the form
$s-t$ where $s \in S$ and $t \in T$. If $x \in [0,n-1]$, $x+S$ will be short for $\{ x \} +S$
and we will sometimes say that $x+S$ is a $\it{rotation}$ of $S$. The
notation $-S$ will be short for the set $\{ 0\}-S$. Given a
graph $G=(V,E)$, the $\it{size}$  of $G$ will equal $|V(G)|$. If $x\in V(G)$, the first neighborhood of $x$ will be the set of vertices $y$ for which $\{ x,y\} \in E$, and
the second neighborhood of $x$ will be the set of $z$ not in the
first neighborhood of $x$ for which there is a $y$ in the first
neighborhood of $x$, with $\{ y,z \} \in E$. The set $N\sb x (G,1)$
will denote the union of $\{ x\}$ and the first neighborhood of $x$, and
the set $N\sb x (G,2)$ will denote the union of $N\sb x (G,1)$ and the
second neighborhood of $x$. If $S \subset V(G)$ the subgraph of $G$
induced by $S$ is the subgraph consisting of the vertex set $S$ and
every edge of $G$ that has both endpoints in $S$.

	Given a positive integer $n \ge 2$ and a set $S \subset [1,\lfloor n/2 \rfloor]$, the circle graph on $n$ vertices   
determined by $S$ will denote the graph $C$ with vertex set $[0,n-1]$,
having two vertices $x,y$ joined by an edge if and only if either
$x-y \in S$ or $y-x \in S$. Note that $N\sb 0(C,2)= \{ 0\} \bigcup S \bigcup (-S)
\bigcup (S+S) \bigcup (-S+S) \bigcup (-S-S)$. From this it follows that
$x\in N\sb 0(C,2)$ if and only if $-x \in N\sb 0 (C,2)$ and $x\in N\sp c \sb
0(C,2)$ if and only if  $-x\in N\sp c \sb 0 (C,2)$. If $C=([0,n-1],E)$ is a circle
graph, any rotation of $[0,n-1]$
defines a graph automorphism of $C$. Therefore, if $H$ is a largest size diameter two subgraph
of $C$, we may assume $0 \in V(H)$ because if it is not a suitable
rotation will show that there is a diameter two subgraph of $C$ with the same
number of vertices as $H$ that contains $0$. It also follows that $N\sb x(C,2)=x+N\sb
0(C,2)$ and since $x+[0,n-1]=[0,n-1]$ that $N\sp c\sb x(C,2)=x+N\sp c\sb 0(C,2)$. Defining $J=J(C)$ to be the circle graph on $n$ vertices determined by $N\sp c \sb 0(C,2) \bigcap  [1,\lfloor n/2
\rfloor]$, we see that $\{ x,y \} \in E(J)$ if and only if there is no $xy$ path in
$C$ containing at most two edges. Thus $V(H)$ is an independent set
in $J(C)$ ( or equivalently $[0,n-1]-V(H)$ is a vertex cover in $J(C)$).
Therefore lower bounds on the size of vertex covers of $J(C)$ imply
upper bounds on the size of diameter two subgraphs of $C$. We rely heavily on the properties of circle graphs throughout the rest of the
paper. In particular, we repeatedly use the fact that if $z \in N\sp c \sb 0 (C,2)$ then for
every $x\in [0,n-1]$, at most one of the two vertices $x$ and $x+z$ can
be in any diameter two subgraph of $C$ because $x+z \in N\sp c \sb x(C,2)$.

	Let $n=2sk+1$ where $s$ is some positive integer and $k \ge 3$ is
the number of colors used. Note that the circle graph determined by $[1,ks]$ is the complete graph on $n$ vertices. The strategy in our constructions will be to partition $[1,ks]$ into $k$ sets $N\sb 1$,\ldots,$N\sb k$ and to prove for every $1 \le i \le k$ that the circle graph
$C\sb i$ determined by $N\sb i$ has no diameter two subgraph on more
than $2s+7$ vertices.

 	Below are a series of Lemma's used in proving the main result when
$k=3$ colors are used.

\begin{lemma} Let $s$ be a positive integer. The set $[2s+1,3s]$ determines a circle graph $C$ on $6s+1$
vertices whose minimum cardinality vertex cover contains at least 
$4s$ vertices. 
\end{lemma}
\underline{Proof of Lemma $2$.} 
Let $D$ be a minimum vertex cover of $C$. By the properties of circle  
graphs we may assume $0 \notin D$. Since $N\sb 0(C,1)=\{ 0\} \bigcup
[2s+1,4s]$, it follows that $[2s+1,4s] \subset D$. Also,  
$\{ x,x+4s\} \in E(C)$ for every $x \in [1,2s]$ and so at least one of
$x$ and $x+4s$ appears in $D$. This
implies $|D| \ge 4s$ and proves Lemma $2$.

\begin{lemma} Let $s$ be a multiple of $3$. The set $[2s/3,s]$ determines a circle graph $C$ on $6s+1$
vertices whose minimum vertex cover contains at least $4s$ vertices.
\end{lemma}

\underline{Proof of Lemma $3$.} 
	Let $t=s/3$. We may assume that $0 \notin D$,
where $D$ is a minimum cardinality vertex cover of $C$. Moreover, by using the
fact that $C$ is a circle graph more carefully, we can by
suitable rotation insure that there is a partition of $[0,n-1]$ into maximal nonempty intervals such that: 
 
 $1$) $[0,n-1]=[e\sb 1,e\sb 2]\bigcup [d\sb 1,d\sb 2]\bigcup \ldots [e\sb
 {2l-1},e\sb {2l}] \bigcup [d\sb {2l-1},d\sb {2l}]$.
  
$2$) $0=e\sb 1 < e\sb 2 +1 =d\sb 1 < d\sb 2+1 =e\sb 3 < \ldots <e\sb {2l} +1
=d\sb {2l-1} \le d\sb {2l} = n-1$.
   
3) $[e\sb {2i-1},e\sb {2i}] \bigcap D= \emptyset$ and $[d\sb
{2l-1},d\sb {2l}] \subset D$ for every $1 \le i \le l$.

	We will denote the length $|[e\sb {2i-1},e\sb {2i}]|$ of the
interval $[e\sb {2i-1},e\sb {2i}]$ by $x\sb i$ and the length of the
interval $[d\sb {2i-1},d\sb {2i}]$ by $y\sb i$.

Now $y\sb i \ge t+1$ for every
$1 \le i \le l$. To see this, let $[d\sb {2i-1},d\sb {2i}]$ be an
arbitrary maximal interval in $D$ and let $w$ be any vertex in this
interval. Because $D$ is a minimum cardinality vertex cover there must
be a vertex $z \in [0,n-1]-D$ such that $ \{ w,z \} \in E(C)$. But $z\in
[0,n-1]-D$ implies that $(z+[2t,3t]) \bigcup (z-[2t,3t]) \subseteq D$.
Certainly $[d\sb {2i-1},d\sb {2i}] \bigcap ( (z+[2t,3t]) \bigcup
(z-[2t,3t]) ) \not= \emptyset$. Since $[d\sb {2i-1},d\sb {2i}]$ is a
maximal interval in $D$, either $z+[2t,3t] \subseteq [d\sb {2i-1},d\sb
{2i}]$ or $z-[2t,3t] \subseteq [d\sb {2i-1},d\sb {2i}]$ whence $y\sb i \ge
t+1$. Moreover, either $y\sb i \le 2t-2$ or $y\sb i \ge 3t$ for every
$1 \le i \le l$. Otherwise, $2t-1 \le y\sb i \le 3t-1$ implies that $2t-2 \le
d\sb {2i}-d\sb {2i-1} \le 3t-2$ which implies that $2t \le (d\sb {2i}+1)
-(d\sb {2i-1}-1) \le 3t$ which is a contradiction since $d\sb {2i-1}-1
$ and $d\sb {2i}+1$ are both in $[0,n-1]-D$. 
\vspace{.2 in}

($1$)	We may assume $y\sb i \le 2t-2$ for every $1 \le i \le l$.
\vspace{.2 in}

\underline{Proof of ($1$).}  First I will show that for every integer $q \in
[0,n-1]$ at least $3t$ vertices of the interval $[q,q+5t-1]$ must be
in $D$. Fix $q$ and let $p,r$ be the least and greatest elements of
$[q,q+2t-1]$  that are in $[0,n-1]-D$. Now $p,r \in [0,n-1]-D$ imply
that $[p+2t,p+3t] \bigcup [r+2t,r+3t] \subset D$. Moreover, for every $1
\le i \le q-p-1$, at least one of $p+i$ and $p+i+3t$ is in $D$. It
follows that every vertex in $[p+2t,r+3t]$ is either in $D$ or
associated with a unique mate in $[p+1,r-1]$ that is in $D$. Moreover
the definition of $p,r$ implies that $[q,p-1] \bigcup [r+1,q+2t-1]
\subset D$. Thus, at least $3t$ vertices of $[q,q+5t-1]$ are in $D$. Now
suppose that some $y\sb i \ge 3t$. Let $[f\sb 1,f\sb 2]$ be a
subinterval of $[d\sb {2i-1},d\sb {2i}]$ with length exactly $3t$. By the
result just proved, the three pairwise disjoint intervals $[f\sb 2+1,f\sb 2+5t]$,$[f\sb 2+5t+1,f\sb 2 + 10t]$ and $[f\sb 2 +10t+1,f\sb 2+15t]$ 
each contribute at least $3t$ vertices to $D$. Since all three are disjoint from
$[f\sb 1,f\sb 2]$, $|D| \ge 12t$ as desired. This proves ($1$). 

	Note that $y\sb i \le 2t-2$ implies $x\sb i +y\sb i + x\sb {i+1} \le
2t$. To see this, note that for every integer $p$ where $y\sb i +1 \le p
\le x\sb i + y\sb i + x\sb {i+1}-1$ there is an $x \in [e\sb {2i-1},e\sb {2i}]$ and
a $w \in [e\sb {2i+1},e\sb {2i+2}]$ such that $p=w-x$. Since $y\sb i \le
2t-2$ for every $1 \le i \le l$, this implies that $2(x\sb 1 +x\sb 2+\ldots
+x\sb l)+(y\sb 1+\ldots+y\sb l) \le
2lt$. Since $y\sb i \ge t+1$ for every $1 \le i \le l$, we may assume $l
\le 11$ because $|D|=y\sb 1 + \ldots +y\sb l$. Also, we have that
$x\sb 1+\ldots +x\sb l \le lt/2$. Since $y\sb 1+\ldots +y\sb l= 18t+1-(x\sb 1+\ldots +x\sb l)$, $|D|=y\sb 1+\ldots +y\sb l \ge (25t/2)+1$ as desired.
This proves Lemma $3$.

\begin{lemma} $[(5s/3)+1,8s/3)]$ determines a circle graph $C$
on $6s+1$ vertices whose largest diameter two subgraph contains at most $2s+3$ vertices.
\end{lemma}
\underline{Proof of Lemma $4$.}  Let $t=s/3$. It can be verified that $N\sp c \sb 0 (C,2)= [8t+1,10t]$.
Now  $V(H) \subset N\sb 0 (C,2)=\{ 0\} \bigcup N \bigcup S \bigcup \bar N \bigcup \bar S$ where $N=[5t+1,8t], S=[13t+1,18t]$, $\bar N =[10t+1,13t]$ and $\bar S=[1,5t]$. We claim that at most $s+1$ vertices of the
$8t$ vertices of $N \bigcup S$ can appear in $H$. To prove this, first assume
that $N \bigcap V(H) = \emptyset$. The set $[13t+1,15t+1]$ is not in the
first neighborhood of $0$ nor is it in the first neighborhood of any
vertex in $\bar N$. Since $0 \in V(H)$, this shows that $|N|+|[13t+1,15t+1]|= 5t+1$ vertices cannot appear in $H$. 
So suppose that $N\bigcap V(H) \not= \emptyset$ and 
let $x$,$y$ be the least and greatest elements of $N$ that are included in
$H$. It suffices to show that at least $5t-1$ of the vertices in $N\bigcup S$
don't appear in $H$. Since $N\sp c \sb 0 (C,2)=[8t+1,10t]$, including $x$ and $y$ in $H$ assures
that $[x+8t+1,x+10t] \bigcap V(H) = \emptyset = [y+8t+1,y+10t]
\bigcap V(H)$.
Also, since $10t \in N\sp c\sb 0(C,2)$, at most
half of the vertices in the set $[x+1,y-2t] \bigcup [x+10t+1,y+8t]$ can
appear in $V(H)$. Because $([5t+1,x-1] \bigcup [y+1,8t])
\bigcap V(H) = \emptyset$, the claim that at most $s+1$ of the $8t$
vertices of $N \bigcup S$ is proved. 

A similar claim holds for $\bar N$ and $\bar S$ whence $|V(H)|
\le 2s+3$. This completes the proof of Lemma $4$.

\begin{theorem} Let $s \ge 9$ be a multiple of $3$ and let $n=6s+1$.
Then the complete graph on $n$ vertices can be
partitioned into $3$ disjoint circle graphs $C\sb 1$, $C\sb 2$ and $C\sb
3$ such that for each $i=1,2,3$, the largest diameter two subgraph of $C\sb
i$ contains at most $2s+3$ vertices.
\end{theorem}

\underline{Proof.}
Let $N\sb 1=[1,s]$, $N\sb 2=[s+1,5s/3] \bigcup [(8s/3)+1,3s]$, and $N\sb 3=[(5
s/3)+1,(8s/3)]$ and let $C\sb i$ denote the circle graph determined by
$N\sb i$ for $1 \le i \le 3$. Let $H\sb i$ be a largest size diameter
two subgraph of $C\sb i$ ($i=1,2,3$) and let $s=3t$.

	It can be verified that $N\sp c \sb 0 (C\sb 1,2)=[2s+1,4s]$ so by Lemma $2$, every vertex cover of $J(C\sb 1)$ contains
at least $4s$ vertices. Since $V(H\sb 1)$ must be an independent set in
$J(C\sb 1)$ and complements of independent sets are vertex covers we
have $|V(H\sb 1)| \le 2s+1$.

	As for $C\sb 2$, note that $N\sb 0 (C\sb 2,1)=
\{ 0\} \bigcup L \bigcup M \bigcup \bar L$ where $L= [3t+1,5t]$, $\bar L=[13t+1,15t]$
and $M= [8t+1,10t]$. Now $L+L=[6t+2,10t]$, $L+M =[11t+2,15t]$, $L+\bar
L=[16t+2,18t] \bigcup [0,2t-1]$, $M+M=[16t+2,2t-1]$, $M+\bar L=[3t+1,7t-1]$ and
$\bar L + \bar L= [8t+1,12t-1]$. From this it can be
verified that $N\sp c \sb 0 (C\sb 2,2)= [2t,3t] \bigcup [15t+1,16t+1]$.
This shows that $J(C\sb 2)$ is the circle graph 
of the hypothesis of Lemma $3$. Applying the conclusion of Lemma $3$ we
get $|V(H\sb 2)| \le 2s+1$.

	Lemma $4$ shows $|V(H\sb 3)| \le 2s+3$ which completes the proof of
Theorem $2$.

\section{ Four or more colors}

\begin{theorem} For all positive integers $k \ge 4$ and $s$, the complete
graph $G$ on $n=2sk+1$ vertices can be decomposed into $k$ circle graphs $C\sb
1$,\ldots,$C\sb k$ each on $n$ vertices such that $E(C\sb 1)$,\ldots,$E(C\sb k)$ forms a partition of the edges of $G$ and such
that for every $1 \le i \le k$, the largest diameter two subgraph of
$C\sb i$ contains at most $2s+7$ vertices.
\end{theorem}
  
\underline{Proof of Theorem $3$.} For $j=1,2$ define $C\sb j$ to be
the circle graph on $2sk+1$ vertices determined by $[j,2s-2+j]\sb 2$.
For integers $j$ satisfying $2 \le j \le k-1$ and $j \notin \{(2k/3)-2,(2k/3)-1\}$ define $C\sb {j+1}$ to be the circle graph on $n$ vertices determined by $[js+1,(j+1)s]$. When $2k/3$ is an
integer define $C\sb {2k/3 -1}$ to be the circle graph on $n$ vertices determined by
$[(2k-6)s/3+2,2ks/3]\sb 2$, and $C\sb {2k/3}$ to be the circle graph on $n$ vertices determined by $[(2k-6)s/3+1,2ks/3-1]\sb 2$. It can be verified
that $E(C\sb 1)$,$E(C\sb 2)$,$\ldots$,$E(C\sb k)$ forms a partition of the
edges of $G$. 
The proof will be broken up into a number of claims. In many cases the proofs
that $C\sb i$ has no diameter two subgraph on more than $2s+7$ vertices
will actually imply the stronger result that 
every vertex cover of $J(C\sb i)$ contains at least $(2k-2)s-6$
vertices.

	First we will prove the Lemma that covers a majority of the
possibilities.

\begin{lemma} Let $k \ge 5$ and $s$ be positive integers and suppose that $j$ is an integer such that $j \notin \{ (2k-3)/3,(2k-2)/3,(2k-1)/3,2k/3
\}$ and $2 \le j \le k-2$. Then the circle graph $D$ on $n=2sk+1$ vertices
determined by
$[js+1,(j+1)s]$ has no diameter two subgraph containing more than $2s+5$
vertices.
\end{lemma}

\underline{Proof of Lemma $5$.} Let $H$ be a largest diameter two
subgraph of $D$. Note that we may assume contains $0$. Let
$A=[1,s-1]$, $B=[js+1,(j+1)s]$, $C=[2js+2,(2j+2)s]$, $\bar C=[(2k-2j-2)s+1,(2k-2j)s-1]$, $\bar B =[(2k-j-1)s+1,(2k-j)s]$ and $\bar A =[(2k-1)s+2,2ks]$. It can be verified that $N\sb 0
(D,2)=\{ 0\} \bigcup A \bigcup B \bigcup C \bigcup \bar C \bigcup \bar B
\bigcup
\bar A$. If $X,Y$ are sets of integers we write $X < Y$ to mean max $X <$ min
$Y$. Let
$B\sp{-}=[(j-1)s+1,js]$ and $B\sp{+}=[(j+1)s+1,(j+2)s]$. If $2 \le j \le (k/2)-1$ then $A < B\sp{-}< B < B\sp {+}< C < \bar C < \bar B < \bar A$. If $j=(k-1)/2$ then $A < B\sp
{-} < B < B\sp {+} < C = \bar C+1 < \bar B < \bar A$. If $k/2 \le j \le (2k-4)/3-1$ then
$A < B\sp {-} < B < B\sp {+} < \bar C < C < \bar B < \bar A$. Finally, if $(2k+1)/3 \le j \le k-2$
then $A < \bar C < B\sp {-} < B < B\sp {+} < \bar B < C < \bar A$. In any case, $B\sp
{-} \bigcup B\sp {+}$ is disjoint from $N\sb 0 (D,2)$. If $V(H) \bigcap
(B \bigcup \bar B) = \emptyset$ then $V(H)=\{ 0 \}$ and we would be done.
From the symmetry we may assume that $B \bigcap V(H) \not= \emptyset$.

	Let $x=js+1+l\sb x \le y=(j+1)s-l\sb y$ ($0 \le l\sb x,l\sb y \le s-1$)
be the least and greatest elements of $V(H) \bigcap [js+1,(j+1)s]$ if they exist and let $\bar y = n-js-s+l\sb {\bar y} \le \bar x=n-js-1-l\sb {\bar x}$
($0 \le l\sb {\bar y},l\sb {\bar x} \le s-1)$ be
the least and greatest elements of $V(H) \bigcap [n-js-s,n-js-1]$
if they exist. Let $t=|[2js+2,2js+2s] \bigcap [n-2js-2s,n-2js-2]|$.
	 
	First it will be shown that we may assume $B
\bigcap V(H) \not= \emptyset \not= \bar B \bigcap V(H)$. 
Suppose that $\bar B \bigcap V(H) = \emptyset$.
The conditions on $j$ imply that
$\bar C $ has no element in $N\sb 0 (D,1)$. Also, no element of $[2ks-2js-2s+1,2ks-2js-1]-[2js+2,(2j+2)s]$ is in the first
neighborhood in $D$ of any
vertex in $[js+1+l\sb x,(j+1)s-l\sb y]$. Therefore
$(\bar C - C) \bigcap V(H) = \emptyset$.
In addition, $B\sp {-} \bigcup B\sp {+} \subset
N\sp c \sb 0 (D,2)$ implies that $([2js+2,(2j+1)s-l\sb y] \bigcup
[(2j+1)s+2+l\sb x,(2j+2)s]) \bigcap V(H) = \emptyset$ and that
$([(2k-1)s+2,2ks-l\sb y] \bigcup [1+l\sb x,s-1]) \bigcap V(H) = \emptyset$. Remembering the definition of $x$ and $y$ and that $l\sb x+l\sb y \le s-1$ we
get that at least $6s-4-t$ vertices in $N\sb 0 (D,2)$ are not in $V(H)$, so $|V(H)| \le 2s+1$ as desired. 

	Therefore, assume that $B \bigcap V(H) \not=
\emptyset \not= \bar B \bigcap V(H)$.
Then the following subsets of $N\sb 0 (D,2)$ are not in $V(H)$
for the reason stated :
		 
	(i) $[js+1,js+l\sb x] \bigcup [(j+1)s+1-l\sb y,(j+1)s] \bigcup
	[n-js-s,n-js-s-1+l\sb {\bar y}] \bigcup [n-js-l\sb {\bar
	x},n-js-1]$ ; definition of $x,y,\bar x,\bar y$.
		  
		  (ii) $[2js+2,(2j+1)s-l\sb y]$ ; $[(j-1)s,js] \subset N\sp c
		  \sb 0 (D,2)$
		  and $y=(j+1)s-l\sb y \in V(H)$.
		   
		   (iii) $[(2j+1)s+2+l\sb x,(2j+2)s]$ ; $[(j+1)s+1,(j+2)s]
		   \subset N\sp c
		   \sb 0 (D,2)$ and $x=js+1+l\sb x \in V(H)$.
			
			(iv) $[n-2js-2s,n-(2j+1)s-2-l\sb {\bar x}] \bigcup
			[n-(2j+1)s+l\sb {\bar
			y},n-2js-2]$; analogous to (ii) and (iii).
			 
			 (v) $[(2k-1)s+2,2ks-$min$\{ l\sb y,l\sb {\bar x} \}]$ ;
			 $\bar x,y \in V(H)$
			 and $[(j-1)s,js] \bigcup [(j+1)s+1,(j+2)s] \subset N\sp c
			 \sb 0 (D,2)$.
			  
			  (vi) $[1+$min$\{ l\sb x,l\sb {\bar y}\} ,s-1]$ ; $\bar y,x
			  \in V(H)$ and
			  $[(j-1)s,js] \bigcup [(j+1)s+1,(j+2)s]$.

	Suppose first that $j=(k-1)/2$. Then
$[2js+2,2js+2s]=([n-2js-2s,n-2js-2]+1)$ so
conditions (ii),(iii) and (iv)
imply that $[(k-1)s+1,ks-$min$\{ 1+l\sb {\bar x},l\sb y \}]
\bigcup [ks+1+$min$\{ l\sb {\bar y},1+l\sb x\} ,(k+1)s]$ is excluded
from $H$.  This and conditions (i), (v) and (vi) imply that at least $4s-4$
vertices are excluded from $H$. Since $|N\sb 0 (D,2)|=6s-1$ it follows
that $V(H) \le 2s+3$ as desired. So assume $j \not= (k-1)/2$. Consider the
quantity $q=|(n-1-2js-s)-(2js+s)|=|2k-4j-2|s$. Then $q$ is always a positive even
multiple of $s$ and $q < (2k-1)s$.

	I claim that $\{ q,q+1\} \bigcap N\sp c \sb 0 (D,2)
\not= \emptyset$. For suppose $q \in N\sb 0 (D,2)$. Then because $q$ is
an even multiple of $s$ it must be that $q=(j+1)s$ or $q=(2j+2)s$ or
$q=(2k-j)s$. Suppose first that $q=(2k-j)s=$max $ \bar B$. If $q+1 \in N\sb 0 (D,2)$ then $q+1=(2k-j-1)s+1$ or $q+1=(2ks-2js-2s)+1$ or $q+1=js+1$, since
$\{ js+1,2ks-2js-2s+1,2ks-js-s+1\}$ is the set of elements of $N\sb 0
(D,2)$ that could possibly equal an even multiple of $s$ plus one. The
case $q+1=(2k-j-1)s+1$ is clearly impossible. The other cases are ruled
out because $\bar B > \bar C$ and $\bar B > B$ for every $2 \le j \le
k-2$ imply $q+1 = 1+$ max $\bar B > $max $\{ 2ks-2js-2s+1,js+1\}$. Thus $q+1 \in N\sp c \sb 0 (D,2)$. 

	If $q=(j+1)s$ then the
restrictions on $j$ in the hypothesis of the Lemma force $j=(2k-3)/5$
which implies that $q=(2k+2)s/5 \not= (6k-4)s/5 = (2ks-2js-2s+1)-1$. By
similar considerations as above this means that $q+1 \in N\sp c \sb 0 (D,2)$ as desired. If $q=(2j+2)s$ then the
restrictions on $j$ force $j=(k-2)/3$ whence $q+1=1+(2k+2)s/3 \not= 2ks-2js-2s+1$ and $q+1 \not= 2ks-js-s+1$. Again, this shows that $q+1 \in N\sp c \sb
0 (D,2)$ and proves the claim that $\{ q,q+1 \} \bigcap N\sp c \sb 0
(D,2) \not= \emptyset$.

	Let $q' =$min $(\{ q,q+1\} \bigcap N\sp c \sb 0 (D,2)\ )$. For every $z
\in [0,2sk]$ at most one of the vertices $z$ and $z+q'$ can appear
in $H$. Therefore,

\vspace{.1 in}

	(vii) "Almost" half of $[(2j+1)s-$min$\{ l\sb y,l\sb {\bar x}\}
	+1,(2j+1)s+$min$\{
	l\sb x,l\sb {\bar y}\}] \bigcup [n-(2j+1)s-$min$\{ l\sb y,l\sb {\bar
	x}\}
	,n-(2j+1)s-1+$min$\{ l\sb {\bar y,x}\}]$ is not in $V(H)$;

\vspace{.1 in}

	In (vii) we say "almost" half because if $q \notin N\sp c \sb 0
(D,2)$ and $q+1 \in N\sp c \sb 0 (D,2)$ then we lose $1$ vertex in the
count.
Thus at least
min$\{l\sb {\bar x},l\sb y\} +$min$\{ l\sb x,l\sb {\bar
y}\} -1$ vertices from the set in (vii) are not in $V(H)$.
Using (i)-(vii) we see that at least $6s-8-t$ vertices are
excluded from $H$ which implies $|V(H)| \le 2s+5$ as desired.
This completes the proof of Lemma $5$.

\begin{lemma} Suppose $s$ and $k \ge 4$ are given positive integers
and let $n=2sk+1$. Then the circle graph $D$ on $n$ vertices
determined by $N=[1,2s-1]\sb 2$ has no diameter two subgraph $H$ on more
than $2s+3$ vertices.
\end{lemma}
 
\underline{Proof of Lemma $6$.}
	Throughout this proof we will abbreviate $[a,b]\sb 2$ by $[a,b]$.
Let $H$ be a largest diameter two subgraph of $D$. Without loss of
generality $H$ contains $0$.
Now $N\sb 0 (D,2)= [(2k-4)s+3,2ks-1] \bigcup [0,4s-2] \bigcup
[1,2s-1] \bigcup [(2k-2)s+2,2ks]$ so $2s+1,(2k-4)s-1,(2k-2)s \in
N\sp c \sb 0 (D,2)$.
Thus at most one third of the vertices from each of the following two
subsets of $N\sb 0 (D,2)$ can appear in $H$:
  
	   $1$) $[1,2s-3] \bigcup ( [1,2s-3]+2s+1) \bigcup ([1,2s-3]+(2k-2)s) $
			   
	   $2$) $[4,2s] \bigcup ( [4,2s]+(2k-4)s-1 ) \bigcup
					 ([4,2s]+(2k-2)s)$
							  
	Therefore at least $4s-4$ vertices are excluded from $H$ 
so $|V(H)| \le |N\sb 0 (D,2)|-(4s-4) = 2s+3$ as desired. This proves
Lemma $6$.

\begin{lemma} Suppose $s$ and $k \ge 4$ are given positive integers
and
that $n=2sk+1$. Then the circle graph $D$ on $n$ vertices
determined by $N=[2,2s]\sb 2$ has no diameter two subgraph $H$ on more
than $2s+1$ vertices.
\end{lemma}
 
\underline{Proof of Lemma $7$.} 
Let $H$ be a largest diameter two subgraph of $D$. We may assume $H$
contains $0$. It is easily verified that $N\sb 0 (D,2)= [(2k-4)s+1,2ks-1]\sb 2
\bigcup [0,4s]\sb 2$. Note that $|N\sb 0 (D,2)|=4s+1$. Since
$(2k-4)s-1 \in N\sp c \sb 0(D,2)$, at most half of the set $[2,4s]\sb 2 \bigcup
[(2k-4)s+1,2ks-1]\sb 2$ is in $V(H)$.  Thus, $H$ contains at most
$2s+1$ vertices. This proves Lemma $7$. 

\begin{lemma} Suppose $s$ and $k \ge 4$ are given positive integers
and that $n=2sk+1$. Then the circle graph $D$ on $n$ vertices 
determined by $N=[(k-1)s+1,ks]$ has no diameter two subgraph on more
than $2s+3$ vertices.
\end{lemma}
 
\underline{Proof of Lemma $8$.}  Let $H$ be a largest monochromatic
diameter two subgraph in $D$. Without loss of generality $H$ contains
$0$. It is easily verified that $N\sb 0 (D,2)= [(2k-2)s+2,2sk] \bigcup
 [0,2s-1] \bigcup [(k-1)s+1,(k+1)s]$.
Since $(k-1)s$ and $2(k-1)s$ are elements of $N\sp c \sb 0 (D,2)$
at most one third of the set
$[2,2s-1] \bigcup ((k-1)s+[2,2s-1]) \bigcup ((2k-2)s+[2,2s-1])$
can appear in $V(H)$. Hence $|V(H)| \le 2s+3$. This proves Lemma $8$.

\begin{lemma} Suppose $k \ge 6$ is an integer divisible by $3$ and $s
$ is a positive integer. Then the circle graph $D=C\sb{2k/3+1}$ on
$2sk+1$ vertices determined by $[(2ks/3)+1,((2k+3)s/3]$ contains no diameter two subgraph on more than $2s+7$ vertices.
\end{lemma}

\underline{Proof of Lemma $9$.} Let $j=2k/3$ and let $H$ be a largest
diameter two subgraph of $D$. Without loss of generality $H$ contains
$0$. Clearly $N\sb 0(D,2) =
[(3j-1)s+2,3js] \bigcup [0,s-1] \bigcup [(j-2)s+1,js-1] \bigcup [js+1,(j+1)s] \bigcup [(2j-1)s+1,2js] \bigcup [2js+2,(2j+2)s]$.
Define for $1 \le i \le s-1$, $v\sb {1,i}=i$, $v\sb {2,i}=(j-2)s+i$, $v\sb {3,i}=js+i$, $v\sb {4,i}=2js+1+i$ and
$\bar v\sb {1,i}=(3j-1)s+1+i$, $\bar v\sb {2,i}=(2j+1)s+1+i$, $\bar v\sb
{3,i}=(2j-1)s+1+i$ and $\bar v\sb {4,i}=(j-1)s+i$. Let $K\sb i$ be the subgraph of $J(D)$
induced by $ \{ v\sb{1,i},v\sb{2,i},v\sb{3,i},v\sb{4,i} \}$ and 
and let $\bar K\sb i$ be the  subgraph of $J(D)$ induced by $\{ \bar
v\sb {1,i},\bar v\sb {2,i},\bar v\sb {3,i},\bar v\sb {4,i} \}$.

It is easy to verify that
$ \{ 2s,(j-2)s,js,(j+2)s+1,2js+1 \} \subset N\sp c \sb 0(D,2)$. 
Therefore, for every $1 \le i \le s-1$, $K\sb i$ ($\bar K\sb i$)
consists of a complete graph on $4$ vertices with the edge between
$v\sb{3,i}$ and $v\sb {4,i}$  ($\bar v\sb {3,i}$ and $\bar v\sb {4,i}$) deleted.For every $1 \le i \le s-1$
let $B\sb i = V(K\sb i)-V(H)$ and $\bar B\sb i = V(\bar K\sb i)-V(H)$.
By definition of $J(D)$, $B \sb i$ ($\bar B\sb i$) is a vertex cover of $K\sb i$ ($\bar K\sb i$).
Note that $|B\sb i|\ge 2$ ($|\bar B\sb i|\ge 2$) with equality if and only if $v\sb{3,i},v\sb{4,i} \in  
V(H)$ ($\bar v\sb{3,i},\bar v\sb{4,i} \in V(H)$). This fact is used
repeatedly in proving the following propositions. 

\vspace{.15 in}
Let $i' \in [1,s-1]$. Then 

\noindent	($1$). 	$v\sb{1,i'} \in V(H)$  implies

		($1$a). $\bar v\sb{1,i} \in \bar B\sb i$ for $1 \le i \le i'$.

		($1$b). $\bar v\sb{3,i} \in \bar B\sb i$ for $1 \le i \le i'-1$.

\noindent	($2$).	$v\sb{2,i'} \in V(H)$ implies 

		($2$a). $v\sb{4,i} \in B\sb i$ for $1 \le i \le s-1$.

		($2$b). $ |B\sb i| \ge 3$ for $1 \le i \le s-1$.

		($2$c). $|B\sb {i'}|+|\bar B\sb{i'}|+\ldots+|B\sb{s-1}|+|\bar
		    B\sb{s-1}| \ge 6(s-i')$.

		($2$d). $v\sb{3,i} \in B\sb i$ for $1\le i \le i'$.

		($2$e). $\bar v\sb{4,i} \in \bar B\sb i$ for $i' \le i \le s-1$.

\noindent ($3$).	$\bar v\sb{2,i'} \in V(H)$ implies

		($3$a). $\bar v\sb{4,i} \in \bar B\sb i$ for $1 \le i \le s-1$.

		($3$b). $|\bar B\sb i| \ge 3$ for $1 \le i \le s-1$.

		($3$c). $|B\sb 1|+|\bar B\sb 1|+\ldots+|B\sb {i'}|+|\bar B\sb {i'}|
		    \ge 6i'$.

		($3$d). $\bar v\sb{3,i} \in \bar B\sb i$ for $i' \le i \le s-1$.

		($3$e). $v\sb{4,i} \in B\sb i$ for $1 \le i \le i'$.

\noindent	($4$).	$v\sb{4,i'} \in V(H)$  implies

		($4$a). $v\sb{3,i'+1} \in B\sb {i'+1}$.

		($4$b). $\bar v\sb{4,i} \in \bar B\sb {i}$ for $1 \le i \le i'$.

\noindent   ($5$).	$\bar v\sb{4,i'} \in V(H)$  implies

		($5$a). $\bar v\sb{3,i'-1} \in B\sb {i'-1}$.
\vspace{.15 in}

\underline{Proof of ($1$)-($5$)} Following each assertion will be the integers in $N\sp c\sb 0(D,2)$ that prove the assertion. ($1$a),($2$d),($2$e),($3$d) and $(3$e),  $[s,2s]$ ;($1$b), $[(2j-2)s+1,(2j-1)s]$ ; ($2$a) and ($3$a), $[(j+1)s+1,(j+3)s]$ ; ($4$a) and
($5$a), $js$ ; ($4$b), $[(j+1)s+1,(j+2)s]$.

	Also, ($2$b) follows from ($2$a). By ($2$b) and ($2$e) we get ($2$c). The
assertion ($3$c) follows similarly. This proves ($1$)-($5$). 
\vspace{.15 in}

($6$)	Suppose that $v\sb {1,i},v\sb {2,i} \in B\sb i$ for every
$i \in [i',i'']$. Then $|B\sb {i'}|+|B\sb {i'+1}|+\ldots+|B\sb {i''}| \ge
3(1+i''-i')-1$. Similarly, if $\bar v\sb {1,i},\bar v\sb{2,i} \in \bar
B\sb i$ for every $i \in [i',i'']$ then $|\bar B\sb {i'}|+\ldots+|\bar B\sb
{i''}| \ge 3(1+i''-i')-1$.
\vspace{.15 in}

\underline{Proof of ($6$).}  The hypothesis implies that $V(H) \bigcap V(K \sb i)
\subset \{ v\sb {3,i},v\sb{4,i} \}$ for $i' \le i \le i''$
Consider all of the indices $p\sb i$,
$i' \le p\sb 1< p\sb 2 <\ldots<p\sb q \le i''$  such that $|B\sb {p\sb i}|=2$ for $i \in [1,q]$. If $q \le 1$ then we are done. Otherwise,
let $r \in [1,q-1]$.
Claim: there is an index $u\sb r$ such that $p\sb r < u\sb r < p\sb {r+1}$
and $|B\sb {u\sb r}|=4$. Suppose not, then ($4$a) implies that 
$V(H) \bigcap V(K\sb {p\sb r+1}) = v\sb {4,p\sb r+1}$ and $V(H) \bigcap
V(K\sb {p\sb {r+1}-1}) = v\sb {3,p\sb {r+1}-1}$. This  implies that
there is an index $u'$ with $p\sb r \le u' < p\sb {r+1}$ such that $v\sb{4,u'},v\sb{3,u'+1} \in V(H)$ which is a contradiction to ($4$a). Then for $r \in [1,q-1]$ we have $|B\sb {p\sb r}|+|B\sb {u\sb r}| =6$
and the result follows. A similar argument proves the other part of   
($6$).

\vspace{.1 in}
($7$)	We may assume that for every $1 \le i \le s-1$, $v\sb{2,i}\in B\sb i$ and $\bar v\sb{2,i} \in \bar B\sb i$.
\vspace{.05 in}

\underline{Proof of ($7$)} 
By ($2$b) and ($3$b) we may assume there is an integer $i' \in [1,s-1]$ such that exactly one of the following two propositions is true: 

a) $v\sb{2,i'} \in V(H)$ and $\{ \bar v\sb{2,1},\ldots,\bar v\sb{2,s-1}\}
\bigcap V(H) = \emptyset$ 

b) $\bar v\sb{2,i'} \in V(H)$ and $\{ v\sb{2,1},\ldots,v\sb{2,s-1} \} \bigcap V(H) = \emptyset$.

Otherwise $|B\sb 1|+|\bar B\sb
1|+\ldots+|\bar B\sb {s-1}| \ge 6s-6$ or $\{v\sb{2,1},\bar
v\sb{2,1},\ldots,\bar v\sb{2,s-1}\} \bigcap V(H)= \emptyset$ and either way
we would be done.

	Assume first that a) is true and let $i'$ be the least integer in
$[1,s-1]$ such that $v\sb{2,i'} \in V(H)$. 
By ($2$a),($2$b) and ($2$c) we have that $v\sb{4,i}
\in B\sb i$ for every $i \in [1,s-1]$,$|B\sb i| \ge 3$ for $i \in [1,s-1]$,
$|B\sb {i'}|+|\bar B\sb {i'}|+\ldots
+|\bar B\sb {s-1}| \ge 6(s-i')$, and hence we may assume that $i'>1$. 
By ($2$d) we have $v\sb{3,i} \in B\sb i$ for
every $1 \le i \le i'$. Thus, for $1 \le i < i'$, $B\sb i = \{
v\sb{1,i},v\sb{2,i},v\sb{3,i},v\sb{4,i} \}$ or $B\sb i=\{ v\sb
{2,i},v\sb{3,i},v\sb{4,i} \}$. For indices $i$ in which the first
possibility holds, $|B\sb i|+|\bar B\sb i| \ge 6$. So we may assume
there is an integer $i'' \in [1,i'-1]$ that is the
greatest index such that the latter holds. By ($1$a) we know
that $\bar v \sb {1,i} \notin V(H)$ for every $1 \le i \le i''$.
From ($6$), this implies that $|B\sb 1|+|\bar B\sb 1|+|B\sb
2|+\ldots+|\bar B\sb {i''}| \ge 6i''-1$ so $|B\sb 1|+|\bar B\sb 1|+\ldots+|\bar
B\sb {s-1}| \ge  6s-7$ whence $|V(H)| \le 2s+4$. A small modification of
this approach works if b) is assumed to hold. This proves ($7$). 

	By ($7$), we may assume that $v\sb{2,i} \in B\sb i$ and $\bar v\sb{2,i}\in \bar B\sb i$ for $1 \le i \le s-1$. Let $l$ denote the greatest
integer in $[1,s-1]$ (if it exists) such that $v\sb{3,l},v\sb{4,l}
\in V(H)$ and let $\bar l$ denote the least integer in $[1,s-1]$ (if
it exists) such that $\bar v\sb{3,\bar l},\bar v\sb{4,\bar l} \in V(H)$. Assume first that they both exist. Note that by ($4$b), $l < \bar l$. Let
$l' < l$ be the greatest integer (if it exists) such that $v\sb{1,l'}
\in V(H)$. By ($1$a),($1$b) and ($4$b), we have $|\bar B\sb i|=4$
for $1 \le i \le l'-1$ and consequently that $|B\sb 1|+|\bar B\sb
1|+\ldots+|B\sb {l'-1}|+|\bar B\sb {l'-1}| \ge 6l'-6$. By ($6$) and
the fact that $l < \bar l$ we have that $|B\sb {l'}|+|\bar B\sb {l'}|+\ldots+|B\sb {l-1}|+|\bar B\sb {l-1}| \ge 6(l-l')-1$ and so $|B\sb 1|+|\bar B\sb
1|+\ldots+|B\sb {l-1}|+|\bar B\sb {l-1}| \ge 6l-7$. If $l'$ did not
exist ($6$) would lead to the same conclusion. Similar reasoning
shows that $|B\sb {\bar l+1}|+|\bar B\sb {\bar l+1}|+\ldots+|B\sb
{s-1}|+|\bar B\sb {s-1}| \ge 6(s-\bar l)-7$. Using the definition of $l$
and $\bar l$, we get that $|B\sb 1|+|\bar B\sb 1|+\ldots+|B\sb {s-1}|+|\bar B\sb{s-1}| \ge 6s-10$ and so $|V(H)| \le 2s+7$. If one or both of $l$ or $\bar l$ does
not exist the above technique can be modified. This proves Lemma $9$.

\begin{lemma} Let $k \ge 6$ be an integer divisible by 3, let $s$ be a positive
integer
 and suppose $j=2k/3$. Then the circle graph
$D=C\sb {2k/3-1}$ on $n=2sk+1$ vertices determined by $[(j-2)s+2,js]\sb 2$ has no diameter two subgraph containing more than $2s+5$ vertices.
\end{lemma}

\underline{Proof of Lemma $10$.}  Let $H$ be a largest diameter two subgraph of $D$. Without loss of generality $H$ contains $0$.  It is easily verified that $N\sb 0 (D,2) = [(3j-2)s+3,2sk-1]\sb 2 \bigcup
[0,2s-2]\sb 2 \bigcup [(j-2)s+2,js]\sb 2 \bigcup [js+1,js+4s-3]\sb 2\bigcup
[(2j-4)s+4,2js]\sb 2 \bigcup [2js+1,(2j+2)s-1]\sb 2$. 
Let $H$ be a largest diameter two subgraph of $D$ containing $0$ and
for an integer $i$ satisfying $1\le i \le s-1$ define $K\sb i$ to be the subgraph of $J(D)$ induced
by $\{ v\sb{1,i},v\sb{2,i},v\sb{3,i},v\sb{4,i} \}$, where   
$v\sb{1,i}=2i$, $v\sb{2,i}=(j-2)s+2i$, $v\sb{3,i}=(2j-4)s+2i+2$, $v\sb{4,i}=js-1+2i$, and let $\bar K\sb i$ be the subgraph of $J(D)$ induced by $\{ \bar
v\sb{1,i}$, $\bar v\sb{2,i}$, $\bar v\sb{3,i}$, $\bar v\sb{4,i} \}$, where $\bar v\sb{1,i}=(3j-2)s+1+2i$, $\bar v\sb{2,i}=2js+1+2i$, $\bar
v\sb{3,i}=(j+2)s-1+2i$ and $\bar v\sb{4,i}=(2j-2)s+2+2i$. Since $\{
2s-1,(j-4)s+3,(j-2)s,js-1,(2j-4)s+2 \} \subset N\sp c \sb 0 (D,2)$ the
graph $K\sb i$ ($\bar K\sb i$) consists of a complete graph on 4
vertices with the edge $\{ v\sb{2,i},v\sb{3,i} \}$ ($\{ \bar
v\sb{2,i},\bar v\sb{3,i} \}$) deleted. For every $1 \le i \le s-1$
let $B\sb i =V(K\sb i)-V(H)$ and $\bar B\sb i = V(\bar K\sb
i)-V(H)$. By definition of $J(D)$, $B \sb i$ ($\bar B\sb i$)
is a vertex cover of $K\sb i$ ($\bar K\sb i$) in $J(D)$. 
Therefore $|B\sb i|\ge 2$ ($|\bar B\sb i|\ge 2$) with equality if and only
if $v\sb{2,i},v\sb{3,i} \in V(H)$ ($\bar v\sb{2,i},\bar v\sb{3,i} \in V(H)$).
 This fact is used repeatedly below.
\vspace{.15 in}

Let $i' \in [1,s-1]$. Then \vspace{.1 in}

\noindent $(1)$. $v\sb{3,i'} \in V(H)$ $(\bar v\sb{3,i'} \in V(H) )$ implies:

 $(1a)$. $\{ v\sb{4,i},\bar v\sb{3,i}\} \bigcap V(H) = \emptyset$ $( 
 \{ v\sb{3,i},\bar v\sb{4,i} \} \bigcap V(H) =\emptyset)$ for $1 \le  i \le s-1$.

 $(1b)$. $\bar v\sb{2,i} \in \bar B\sb i$ for  $1 \le i
 \le i'$. $(v\sb{2,i} \in B\sb i$ for $i' \le i \le s-1)$.

 $(1c)$. $v\sb{2,i'+1} \in B\sb {i'+1}$ $(\bar v\sb{2,i'-1} \in \bar B\sb
  {i'-1})$.

  $(1d)$. $\bar v\sb{4,i'} \in \bar B\sb{i'}$  $(v\sb{4,i'} \in
B\sb{i'})$.

\noindent $(2)$. $v\sb{2,i'} \in V(H) $ $(\bar v\sb{2,i'} \in V(H) )$ implies:

$(2a)$. $\bar v\sb{1,i'-1} \in \bar B\sb {i'-1}$ $(v\sb {1,i'+1} \in B\sb
i)$.

$(2b)$. $\bar v\sb{4,i'} \in \bar B\sb {i'}$ $(v\sb{4,i'} \in B\sb{i'})$.

\noindent $(3)$. $\bar v\sb{1,i'} \in V(H)$ $(v\sb {1,i'} \in V(H) )$ implies 

$(3a)$. $\bar v\sb{4,i} \in \bar B\sb {i}$ for $i' \le i \le s-1$
$(v\sb{4,i}\in B\sb {i}$ for every $1 \le i \le i')$.

\vspace{.15 in}

\underline{Proof of $(1)$-$(3)$.} Following each assertion will be the integers in $N\sp
c\sb 0(D,2)$ that prove the assertion. ($1$a)  $[1,js-1]\sb
2$; ($1$b) $[1,js-1]\sb 2$ ;
($1$c) $(j-2)s$; ($1$d) $2s$ ; ($2$a) $2js-1$ ; ($2$b) $js+2$ ; $(3a)$ 
$[1,js-1]\sb 2$. This proves $(1)$-$(3)$. 

	Consider the two propositions 
 i) There is an $i \in [1,s-1]$ such that$\{ v\sb{2,i},v\sb{3,i} \} \subset V(H)$
and  ii) There is an $i \in [1,s-1]$ such that $\{ \bar v\sb{2,i},\bar v\sb{3,i}\} \subset V(H)$. 
By ($1$a) at most 1 of the propositions is true. Also, we may assume at
least 1 is true, otherwise $|B\sb i|+|\bar B\sb i| \ge 6$ for every $i\in
[1,s-1]$ and we would be done. Suppose first that
proposition i) is true and let $1 \le p\sb 1 < p\sb 2 < \ldots < p\sb r \le
s-1$ be the indices for which $\{ v\sb{2,p\sb i},v\sb{3,p\sb i}\}
\subset V(H)$. Clearly $|B\sb i| \ge 3$ for every $i \in [1,s-1]-\{ p\sb
1,\ldots,p\sb r\}$ and since proposition ii) is false, $|\bar B\sb i| \ge 3$ for
every $i \in [1,s-1]$. By ($1$a), ($1$b) and ($1$d), we must have $V(\bar K\sb {p\sb i})
\bigcap V(H) \subset \{ \bar v\sb{1,p\sb i} \}$ for $i \in [1,r]$. If $\bar
v\sb{1,i} \notin V(H)$ 
then $|B\sb {p\sb i}|+|\bar B\sb {p\sb i}| \ge
6$. So let  $i' \in [1,r]$ be the least subindex such that $\bar v\sb{1,p\sb {i'}}
\in V(H)$.  By ($3$a), $\bar v\sb {4,i''} \in B\sb {i''}$ for all
$p\sb {i'} \le i'' \le p\sb r$. This and $1$a), $1$b), $2$a) imply that $|\bar B\sb {p\sb
i-1}| = 4$ for $i \in [i'+1,r]$. By ($1$c) $p\sb {i+1} > p\sb i +1$ for
every $i \in [1,r-1]$ and hence $|B\sb {p\sb i-1}| \ge 3$ for all $i \in [1,r]$. Therefore, for every $i\in [i'+1,r]$, we
have $|B\sb {p\sb i-1}|+|\bar B\sb{p\sb i-1}|+|B\sb {p\sb i}|+|\bar B\sb
{p\sb i}| \ge 12$. This implies that $|B\sb 1|+|\bar B\sb 1|+\ldots+|\bar
B\sb {s-1}| \ge 6s-7$ whence $|V(H)| \le 2s+4$.   
	If proposition ii) is true this technique can be modified slightly
by letting $i'$ be the greatest subindex for which $v\sb{1,p\sb {i'}} \in V(H)$.This proves Lemma $10$.

\begin{lemma} Let $k \ge 6$ be an integer divisible by $3$, let $s$ be a
positive integer
and suppose $j=2k/3$. Then the circle graph
$D=C\sb{(2k/3)}$ on $2sk+1$ vertices determined by $[(j-2)s+1,js-1]\sb 2$ has no diameter two subgraph containing more than $2s+5$ vertices.
\end{lemma}

\underline{Proof of Lemma $11$.}  Let $H$ be a largest size diameter two subgraph of
$D$. Without loss of generality, $H$ contains $0$. It is easily verified that $N\sb 0 (D,2) = [(3j-2)s+3,2sk-1]\sb 2
\bigcup [0,2s-2]\sb 2 \bigcup [(j-2)s+1,js-1]\sb 2 \bigcup [js+3,(j+4)s-1]\sb
2 \bigcup
[(2j-4)s+2,2js-2]\sb 2 \bigcup [2js+2,(2j+2)s]\sb 2$. For $1 \le i \le s-1$ define
$v\sb{1,i}=2i$,
$v\sb{2,i}=(j-2)s-1+2i$,$v\sb{3,i}=js+1+2i$,$v\sb{4,i}=(2j-4)s+2i$,$\bar
v\sb{1,i}=(3j-2)s+1+2i$,$\bar v\sb{2,i}=2js+2i+2$,$\bar
v\sb{3,i}=(2j-2)s+2i$, and $\bar v\sb{4,i}=(j+2)s+1+2i$. For the same
set of $i$ define $K\sb i$ ($\bar K\sb i$) to be the subgraph of $J=J(D)$ induced by the vertex set
$\{ v\sb{1,i},v\sb{2,i},v\sb{3,i},v\sb{4,i} \}$ ($\{ \bar v\sb{1,i},\bar
v\sb{2,i},\bar v\sb{3,i},\bar v\sb{4,i} \}$). Because $\{
2s+2,(j-4)s-1,(j-2)s-1,js+1,(2j-4)s\} \subset N\sp c \sb 0(D,2)$, 
$K\sb i$ ($\bar K\sb i$) consists of a complete graph on 4 vertices with
the edge $\{ v\sb{2,i},v\sb{4,i} \}$ ($\{ \bar v\sb{2,i},\bar v\sb{4,i}
\}$) deleted. For every $i \in [1,s-1]\sb 1$, let $B\sb i = V(K\sb i)-V(H)$ and $\bar B\sb i = V(\bar K\sb i)-V(H)$. By the nature of $J(D)$,  
$B\sb i$ ($\bar B\sb i$) must be a vertex cover of $K\sb i$
($\bar K\sb i$) in $J$.  Note that $|B\sb i|\ge 2$ ($|\bar B\sb i|\ge 2$) with equality if and only if $v\sb{2,i},v\sb{4,i} \in V(H)$ ($\bar v\sb{2,i},\bar v\sb{4,i} \in V(H)$)

 	There can be at most one index $i \in [1,s-1]$ such that $v\sb{2,i},v\sb{4,i}
\in V(H)$ because $[1,(j-2)s-1]\sb 2 \subset N\sp c \sb 0(D,2)$
and $v\sb {2,i} \in V(H)$ implies that $v\sb{4,i'} \in B\sb {i'}$ for
every $i' \in [1,i-1]$. Similar reasoning implies that $\bar
v\sb{2,i},\bar v\sb{4,i} \in V(H)$ for at most one $i \in [1,s-1]$.
Therefore $|B\sb 1|+|\bar B\sb i|+|B\sb 2|+\ldots+|\bar B\sb {s-1}| \ge
6s-8$ which implies $|V(H)| \le 2s+5$, as desired. This proves Lemma $11$.

\begin{lemma} Let $k \ge 4$ and $s$ be a positive integers. If $j\in \{ (2k-2)/3 ,(2k-1)/3 \}$ and $j$ is an integer
then the circle graph $D$ on $n=2sk+1$ vertices determined by $[js+1,(j+1)s]$ has no diameter
two subgraph on more than $2s+7$ vertices.
\end{lemma}

\underline{Proof of Lemma $12$.} Suppose first that $j=(2k-2)/3$. As usual, let $H$
be a largest diameter 2 subgraph of $D$ containing $0$. It is easily
verified that $N\sb 0 (D,2) = [(3j+1)s+2,(3j+2)s] \bigcup [0,s-1] \bigcup
[js+1,(j+2)s-1] \bigcup [2js+2,(2j+2)s]$. Now $js,2js \in N\sp c \sb 0
(D,2)$, so at most one third of each of the following two subsets of $N\sb
0 (D,2)$ can occur in $H$:

1) $[2,s-1] \bigcup (js+[2,s-1]) \bigcup (2js+[2,s-1])$

2) $[(3j+1)s+2,(3j+2)s-1] \bigcup ([(3j+1)s+2,(3j+2)s-1]-js) \bigcup ([(3j+1)s+2,(3j+2)s-1]-2js)$. 

	Thus at least $4s-8$ vertices are excluded from $H$ so $|V(H)| \le
2s+5$.

	Now suppose $j=(2k-1)/3$. Let $H$ be a largest diameter two subgraph
containing $0$. Since $N\sb 0 (D,2) = [3js+2,(3j+1)s] \bigcup [0,s-1] \bigcup
[(j-1)s+1,(j+1)s] \bigcup [2js+1,(2j+2)s]$, $(j-1)s-1,(j+1)s+1,2js
\in N\sp c \sb 0 (D,2)$. Therefore at most one third of each of the
following two subsets of $N\sb 0 (D,2)$ can appear in $H$:

1)  $[2,s-1] \bigcup ( [2,s-1]+(j-1)s-1 ) \bigcup ([2,s-1]+2js)$.

2)  $[3js+2,(3j+1)s-1] \bigcup ( [3js+2,(3j+1)s-1]-(j-1)s+1 ) \bigcup
([3js+2,(3j+1)s-1]-2js)$.

	 Thus $|V(H)| \le 2s+7$ as desired. This proves Lemma $12$.

	Lemma's $5$-$12$ complete the proof of Theorem $3$.  Q.E.D. 

\vspace{.07 in}

I would like to acknowledge the God of the Old and New Testaments because He is my creator, sustainer and redeemer.

\end{document}